\documentclass[11pt,a4paper]{article}
\usepackage[top=3cm, bottom=3cm, left=3cm, right=3cm]{geometry}
\usepackage[latin1]{inputenc}
\usepackage{amsmath}
\usepackage{amsthm}
\usepackage{amsfonts}
\usepackage{amssymb}
\usepackage{graphicx}

\usepackage{amsmath,amsthm,amssymb,graphicx, multicol, array}
\usepackage{enumerate}
\usepackage{enumitem}
\usepackage{hyperref}

\DeclareMathOperator{\ord}{ord}

\DeclareMathOperator{\lcm}{lcm}
\newtheorem{thm}{Theorem}[section]
\newtheorem{lem}{Lemma}[section]

\newtheorem{exa}{Example}[section]

\newtheorem{dfn}{Definition}[section]

\newcommand{\N}{\mathbb{N}}
\newcommand{\Z}{\mathbb{Z}}

\newcommand{\R}{\mathbb{R}}
\newcommand{\C}{\mathbb{C}}

\newcommand{\tP}{\mathbb{P}}

\usepackage{fancyhdr}
\title{The Generalized Artin Primitive Root Conjecture}
\date{}
\author{N. A. Carella}
\begin{document}

\maketitle

\begin{abstract}
Asymptotic formulas for the number of integers with the primitive root 2, and the generalized Artin conjecture for multiplicative subsets of composite integers with fixed admissible primitive roots \(u\neq \pm 1,v^2\), are presented here.
\end{abstract}
\tableofcontents
\vskip .15 in

\text{AMS MSC:} Primary 11A07, Secondary 11Y16, 11M26.\\
\text{Keywords:} Primitive root, Artin primitive root conjecture, Generalized primitive root.

\section{Introduction} \label{s1}
A generalized Artin conjecture for multiplicative subsets of composite integers with fixed primitive roots is presented here. The focus is on developing an asymptotic formula for the number of integers in the multiplicative subset $\mathcal{N}_2$ with the primitive root 2, modulo the generalized Riemann hypothesis, see Definition \ref{dfn1.05}. The analysis easily extends to all the admissible primitive roots \(u\neq \pm 1,v^2\). This analysis spawn new questions about the structure of an \textit{L}-series associated with the multiplicative subset of integers with a fixed primitive root $u$, and related ideas. 

\subsection{Subset of Integers With Fixed Primitive Root 2}
The symbols $\N$, $\tP$, $\Z$, $\R$, and $\C$ denotes the set of natural numbers, the set of primes, the set of integers, the set of real numbers, and the set of complex numbers.  
\begin{dfn} \label{dfn1.01} Let $n \geq 1$ be an integer. The order $\ord_n(u)$ of the element $u\in \Z/n\Z$ is defined by
\begin{equation}
\ord_n(u)=\min \left\{ d: u^d \equiv 1 \bmod n \right\} .
\end{equation} 
\end{dfn}

\begin{dfn} \label{dfn1.03} Let $u\ne \pm 1, v^2$ be an integer. The subset $\mathcal{P}_u$ of primes  is defined by
\begin{equation}
\mathcal{P}_u=\left\{ p\in \mathbb{P}:\ord_p(u)=p-1 \right\} \subset \mathbb{P}.
\end{equation} 
\end{dfn}

\begin{dfn} \label{dfn1.05} Let $u\ne \pm 1, v^2$ be an integer. The multiplicative subset $\mathcal{N}_u$ of integers generated by the subset of primes $ \mathcal{P}_u$ is defined by
\begin{equation} \label{eq2-40}
\mathcal{N}_u =\left\{ n\in \mathbb{N}:\ord_n(u)=\lambda(n) \text{ and } p \mid n \Rightarrow \ord_p(u)=p-1, \ord_{p^2}(u)=p(p-1) \right\}.
\end{equation} 
\end{dfn}

The Artin primitive root conjecture states that the integer 2 is a primitive root mod \(p\) for infinitely many primes. Id est,
\begin{eqnarray}
\mathcal{P}_2&=&\#\left\{ p\leq x:\text{ord}_p(2)=p-1 \right\}=\alpha _2\pi (x) \\
&=&\{ 3, 5, 11, 13, 19, 29, 37, 53, 59, 61, 67, 83, 101, 107, 107,131,139,\text{...} \} \nonumber.
\end{eqnarray} 
Conditional on the generalized Riemann hypothesis, Hooley proved that the subset of primes $\mathcal{P}_2$ has nonzero density \(\alpha _2=\delta \left(\mathcal{P}_2\right)>0\), see Theorem \ref{thm7.1} or \cite{HC67}. Moreover, it has the counting function 
\begin{equation}
\pi _2(x)=\#\left\{ p\leq x:\text{ord}_p(2)=p-1 \right\}=\alpha _2\pi (x).
\end{equation}
Here \(\pi(x)=\#\left\{ p\leq x \right\} \) is the primes counting function. Partial unconditional results on the Artin primitive root conjecture are also available in \cite{GM84}, et alii. The subset \(\mathcal{P}_2\) is utilized here to generate the multiplicative subset of composite integers 
\begin{eqnarray} \label{eq2-20}
\mathcal{N}_2 &=&\left\{ n\in \mathbb{N}:\ord_n(2)=\lambda(n) \text{ and } p \mid n \Rightarrow \ord_p(2)=p-1, \ord_{p^2}(2)=p(p-1) \right\} \nonumber \\
&=&\left\{ 3, 5, 3^2,11, 3\cdot 5,19,5^2,3^3, 29,3\cdot 11,37, 45, 53,55,3\cdot 19,61, \ldots \right\},
\end{eqnarray}
which have \(u=2\) as a primitive root. The constraint $\ord_{p^2}(2)=p(p-1)$ sieves the subset of composite integers generated by the subset of Wieferich primes 
	\begin{equation}
	\mathcal{W}_2=\left\{ p\in \mathbb{P}:2^{p-1} \equiv 0 \bmod p^2 \right\}=\{1093, 3511, \ldots,\}.
	\end{equation}
The underlining structures of the asymptotic counting formulas 
\begin{equation}
N_2(x)=\#\left\{ n\leq x:n \in \mathcal{N}_2 \right\}
\end{equation}
and 
\begin{equation}
N_u(x)=\#\left\{ n\leq x:n \in \mathcal{N}_u \right\}
\end{equation}
for the multiplicative subsets of integers $\mathcal{N}_2$ and $\mathcal{N}_u$ will be demonstrated here. These resulte are consistent with the heuristic explained in \cite[p. 10]{LP02}, and have the expected asymptotic orders \(N_2(x)=o(x)\) and \(N_u(x)=o(x)\). \\

\begin{thm} \label{thm1.1} 
Assuming the generalized Riemann hypothesis, the integer $2$ is a primitive root mod \(n\) for infinitely many composite integers
	\(n\geq 1\). Moreover, the number of integers \(n\leq x\) in $\mathcal{N}_2$ such that $2$ is a primitive root mod \(n\) has the asymptotic formula
	
	\begin{equation}
	N_2(x)=\left(\frac{e^{\gamma_2-\gamma \alpha_2}}{\Gamma \left(\alpha _2\right)}+o(1)\right)\frac{x}{(\log  x)^{1-\alpha _2}}\prod _{p\in \mathcal{W}
	} \left(1-\frac{1}{p^2}\right),
	\end{equation}
where \(\alpha_2>0\) is Artin constant, $\gamma_2$ is a generalized Euler constant, $\Gamma (s)=\int _0^{\infty }t^{s-1}e^{-s t}d t$, where $s\in \mathbb{C}$ is a complex number, is the gamma function, for all large numbers \(x\geq 1\).
\end{thm}

The average order of the counting function \(N_2(x)\) has the same average order as the counting function \(L(x)=\#\{ n\leq x:\mu (\lambda(n))=\pm 1 \}=(\kappa +o(1))x(\log  x)^{\alpha _2-1}\) for the number of squarefree values of the Carmichael function \(\lambda\), but it has a different
density \(\kappa \neq e^{\gamma_2-\gamma \alpha_2} /\Gamma \left(\alpha _2\right)\), see \cite{PS03}. This should be compared to the counting function \(\pi _2(x)=\#\left\{
p\leq x:\ord_p(2)=p-1 \right\}=\alpha _2\pi (x)\) for the set of primes having 2 as a primitive root, and the counting function \(T(x)=\#\{
p\leq x:\mu (\varphi (n))=\pm 1 \}=\alpha _2\pi (x)\) of squarefree values of the Euler totient function \(\varphi\). All these asymptotic formulae are closely related and scaled by the constant \(\alpha _2\). \\

The general asymptotic formula for the number of integers \(n\leq x\) in the multiplicative subset of composite integers $\mathcal{N}_u $ such that $u \ne \pm 1, v^2$ is a primitive root mod \(n\) has the form stated here.

\begin{thm} \label{thm1.2} 
	Assuming the generalized Riemann hypothesis, the integer $u\ne \pm 1, v^2$ is a primitive root mod \(n\) for infinitely many composite integers
	\(n\geq 1\). Moreover, the number of integers \(n\leq x\) in $\mathcal{N}_u$ such that $u$ is a primitive root mod \(n\) has the asymptotic formula
	
	\begin{equation}
	N_u(x)=\left(\frac{e^{\gamma_u-\gamma \alpha_u}}{\Gamma \left(\alpha _2\right)}+o(1)\right)\frac{x}{(\log  x)^{1-\alpha _u}}\prod _{p\in \mathcal{W}
	} \left(1-\frac{1}{p^2}\right),
	\end{equation}
	
	where \(\alpha_u>0\) is Artin constant, $\gamma_u$ is a generalized Euler constant, $\Gamma (s)=\int _0^{\infty }t^{s-1}e^{-s t}d t$, where $s\in \mathbb{C}$ is a complex number, is the gamma function,  and 
\begin{equation}
\mathcal{W}=\left\{ p\in \mathbb{P}:u^{p-1} \equiv 0 \bmod p^2 \right\}
\end{equation}
is the set of Abel-Wieferich primes, for all large numbers \(x\geq 1\)
\end{thm}

The generalized Euler constant, and Mertens constant are discussed in Section \ref{s5}. 
Sections \ref{s2} to \ref{s6} provide some essential background results. The proofs of Theorem \ref{thm1.1} and Theorem \ref{thm1.2} ar settled in Section \ref{s7}. \\

\section{Some Arithmetic Functions} \label{s2}
The Euler totient function counts the number of relatively prime integers \(\varphi (n)=\#\{ k:\gcd (k,n)=1 \}\). This counting function is compactly expressed by
the analytic formula \(\varphi (n)=n\prod_{p \mid n}(1-1/p),n\in \mathbb{N} .\)\\

\begin{lem} {\normalfont (Fermat-Euler)} \label{lem2.1}If \(a\in \mathbb{Z}\) is an integer such that \(\gcd (a,n)=1,\) then \(a^{\varphi (n)}\equiv
	1 \bmod n\).
\end{lem}

The Carmichael function is basically a refinement of the Euler totient function to the finite ring \(\mathbb{Z}/n \mathbb{Z}\). 

\begin{dfn} Given an integer
\(n=p_1^{v_1}p_2^{v_2}\cdots  p_t^{v_t}\), the Carmichael function is defined by
\begin{equation}
	\lambda (n)=\text{lcm}\left (\lambda \left(p_1^{v_1}\right),\lambda \left (p_2^{v_2}\right ) \cdots  \lambda \left (p_t^{v_t}\right ) \right )
	=\prod _{p^v  \mid \mid  \lambda (n)} p^v,
\end{equation}
where the symbol \(p^v \mid \mid n,\nu \geq 0\), denotes the maximal prime power divisor of \(n\geq 1\), and 
\begin{equation}
	\lambda
	\left(p^v\right)= \left \{
	\begin{array}{ll}
		\begin{array}{ll}
			\varphi \left(p^v\right) 
			& \text{ if } p\geq 3\text{ or } v\leq 2,  \\
			2^{v-2} 
			& \text{ if } p=2 \text{ and }v \geq 3. \\
		\end{array}
		
	\end{array} \right . 
\end{equation}
\end{dfn} 
The two functions coincide, that is, \(\varphi(n)=\lambda (n)\) if \(n=2,4,p^m,\text{ or } 2p^m,m\geq 1\). And \(\varphi \left(2^m\right)=2\lambda
\left(2^m\right)\). In a few other cases, there are some simple relationships between \(\varphi (n) \text{ and } \lambda (n)\). In fact, it seamlessly
improves the Fermat-Euler Theorem: The improvement provides the least exponent \(\lambda (n) \mid  \varphi (n)\) such that \(a^{\lambda (n)}\equiv
1 \bmod n\). \\

\begin{lem} \label{lem2.2} { \normalfont (\cite{CR10})}  Let \(n\in \mathbb{N}\) be any given integer. Then
\begin{enumerate} [font=\normalfont, label=(\roman*)]
\item The congruence \(a^{\lambda (n)}\equiv
1 \bmod n\) is satisfied by every integer \(a\geq 1\) relatively prime to \(n\), that is \(\gcd (a,n)=1\).

\item In every congruence \(x^{\lambda (n)}\equiv 1 \bmod n\), a solution \(x=u\) exists which is a primitive root \(\bmod  n\), and for any such
solution \(u\), there are \(\varphi (\lambda (n))\) primitive roots congruent to powers of \(u\).
\end{enumerate} 
\end{lem}

\begin{proof} (i) The number \(\lambda (n)\) is a multiple of every \(\lambda \left(p^v\right)=\varphi\left(p^v\right) \) such that \(p^v \mid  n\). Ergo, for any relatively prime integer \(a\geq 2\), the system of congruences 
\begin{equation}
a^{\lambda (n)}\equiv 1\bmod p_1^{v_1}, \quad a^{\lambda (n)}\equiv 1\bmod p_2^{v_2}, \quad \ldots, \quad a^{\lambda (n)}\equiv 1\bmod p_t^{v_t},
\end{equation}
where \(t=\omega (n)\) is the number of prime divisors in \(n\), is valid. 
\end{proof}

\begin{dfn}
	An integer \(u\in \mathbb{Z}\) is called a \textit{primitive root} \(\text{mod } n \) if the least exponent \(\min  \left\{ m\in \mathbb{N}:u^m\equiv 1 \bmod n \right\}=\lambda
	(n)\).
\end{dfn}

\begin{lem} \label{lem2.4}  {\normalfont (Primitive root test)} An integer $u \in \Z$ is a primitive root modulo an integer $n \in \N$ if and only if 
\begin{equation}\label{eq2-52}
u^{\lambda(n)/p} -1\not \equiv 0 \mod  n
\end{equation}
for all prime divisors $p \mid \lambda(n)$.
\end{lem}
The primitive root test is a special case of the Lucas primality test, introduced in \cite[p.\ 302]{ LE78}. A more recent version appears in \cite[Theorem 4.1.1]{CP05}, and similar sources. 

\begin{lem} \label{lem2.3}  Let $n$, and $u \in \mathbb{N}$ be integers, $\gcd(u,n)=1$. If $u$ is a primitive root modulo $p^k$ for each prime power divisor $p^k \mid n$, then, the integer $u \ne \pm 1, v^2$ is a primitive root modulo $n$.
\end{lem}

\begin{proof} Without loss in generality, let $n=pq$ with $p\geq 2$ and $q\geq 2$ primes. Let $u$ be a primitive root modulo $p$ and modulo $q$ respectively. Then
\begin{equation} \label{eq2-01}
u^{(p-1)/r}  -1\not \equiv 0 \bmod  p  \qquad  \text{  and  } \qquad u^{(q-1)/s} -1 \not \equiv 0 \bmod  q,
\end{equation}
for every prime $r \mid p-1$, and every prime $s \mid q-1$ respectively, see Lemma \ref{lem2.4}. Now, suppose that $u$ is not a primitive root modulo $n$. In particular,
\begin{equation}\label{eq2-02}
u^{\lambda(n)/t} -1\equiv 0 \bmod  n  
\end{equation}
for some prime divisor $t \mid \lambda(n)$.\\

Let $v_t(\lambda(n))$, $v_t(p-1)$, and $v_t(q-1)$ be the $t$-adic valuations of these integers. Since $\lambda(n)=\lcm(\varphi(p-1), \varphi(q-1))$, it follows that at least one of the relations
\begin{equation}\label{eq2-23}
v_t(\lambda(n))=v_t(p-1)  \qquad  \text{ or } \qquad v_t(\lambda(n))=v_t(q-1) 
\end{equation}
is valid. As consequence, at least one of the congruence equations
\begin{equation}\label{eq2-06}
u^{\lambda(n)/t} -1\equiv 0 \bmod  n  \qquad  \Longleftrightarrow \qquad u^{\lambda(n)/t}-1 \equiv \mod p
\end{equation}
or
\begin{equation}\label{eq2-08}
u^{\lambda(n)/t} -1\equiv 0 \bmod  n  \qquad  \Longleftrightarrow \qquad u^{\lambda(n)/t}-1 \equiv \mod q
\end{equation}
fails. But, this in turns, contradicts the relations in (\ref{eq2-01}) that $u$ is a primitive root modulo both $p$ and $q$. Therefore, $u$ is a primitive root modulo $n$. 
\end{proof}

\begin{exa} {\normalfont The integer $2$ is a primitive root molulo both $p=37$ and $q=61$. Let $n=p\cdot q=37 \cdot 61$, $\varphi(p-1)=2^2\cdot 3^2$, $\varphi(q-1)=2^2\cdot 3 \cdot 5$, and $\lambda(n)=\lcm(\varphi(p-1), \varphi(q-1))=2^2\cdot 3^2 \cdot 5$. The corresponding congruences and $t$-adic valuations of these integers are these.
\begin{itemize}
\item For $t=2$, the valuations are: $v_2(\lambda(n))=v_2(p-1)=v_2(q-1)=2$. The assumption that 2 is not a primitive root modulo $n$ is not valid: 
\begin{equation} \label{eq2-72}
u^{\lambda(n)/2} -1\equiv 0 \bmod  n  
\end{equation}
fails because at least one 
\begin{equation} \label{eq2-76}
2^{\lambda(n)/2}-1 \not \equiv 0\mod p  \qquad  \text{ or } \qquad 2^{\lambda(n)/2}-1 \not \equiv 0\mod q
\end{equation}
contradicts it.
\item For $t=3$, the valuations are: $v_3(\lambda(n))=v_3(p-1)=2$, and $v_3(q-1)=1$. The assumption that 2 is not a primitive root modulo $n$ is not valid: 
\begin{equation}\label{eq2-78}
u^{\lambda(n)/3} -1\equiv 0 \bmod  n  
\end{equation}
fails because at least one 
\begin{equation}\label{eq2-79}
2^{\lambda(n)/3}-1 \not \equiv 0\mod p  \qquad  \text{ or } \qquad 2^{\lambda(n)/3}-1 \equiv 0 \mod q
\end{equation}
contradicts it.
\item For $t=5$, the valuations are: $v_5(\lambda(n))=v_5(q-1)=1$, and $v_5(p-1)=0$. The assumption that 2 is not a primitive root modulo $n$ is not valid: 
\begin{equation}\label{eq2-78}
u^{\lambda(n)/5} -1\equiv 0 \bmod  n  
\end{equation}
fails because at least one 
\begin{equation}\label{eq2-79}
2^{\lambda(n)/5}-1  \equiv 0\mod p  \qquad  \text{ or } \qquad 2^{\lambda(n)/5}-1 \not \equiv 0 \mod q
\end{equation}
contradicts it.
\end{itemize}
Since the congruence (\ref{eq2-72}) fails for every prime divisor $t=2,3,5$ of $\lambda(n)=2^2\cdot 3^2 \cdot 5$, it implies that $2$ is primitive root modulo $n=37 \cdot 61$.
}
\end{exa}

\section{Characteristic Function In Finite Rings} \label{s3}
The symbol \(\ord_{p^k}(u)\) denotes the order of an element \(u\in \left(\mathbb{Z}\left/p^k\right.\mathbb{Z}\right)^{\times }\) in the multiplicative
group of the integers modulo \(p^k\). The order satisfies the divisibility condition \(\text{ord}_{p^k}(u) \mid \lambda (n)\), and
primitive roots have maximal orders \(\ord_{p^k}(u)=\lambda (n)\). The basic properties of primitive root are explicated in \cite{AT76}, \cite{RH94}, et cetera. The characteristic function \(f:\mathbb{N}\longrightarrow \{ 0, 1 \}\) of a fixed primitive root \(u\) in the finite ring \(\mathbb{Z}\left/p^k\right.\mathbb{Z}\),
the integers modulo \(p^k\), is determined here. 
\\

\begin{lem} \label{lem3.1}Let \(p^k,k\geq 1\), be a prime power, and let \(u\in \mathbb{Z}\) be an integer such that \(\gcd \left(u,p^k\right)=1\). Then
\begin{enumerate} [font=\normalfont, label=(\roman*)]
\item The characteristic \(f\) function of the primitive root \(u \bmod p^k\) is given by  
\begin{equation}
		f\left(p^k\right) =\left \{
		\begin{array}{ll}
			\begin{array}{ll}
				1 & \text{    }\text{ if } p^k=2^k, k\leq 2, \\
				0 & \text{    }\text{ if } p^k=2^k, k>2, \\
				1 & \text{    }\text{ if } \ord_{p^k}(u)=p^{k-1}(p-1), p>2, \text{for} \text{ any } k\geq 1, \\
				0 & \text{    }\text{ if } \ord_{p^k}(u)\neq p^{k-1}(p-1), \text{ and } p>2, k\geq 1. \\
			\end{array}
		\end{array} \right .
	\end{equation}
\item The function \(f\) is multiplicative, but not completely multiplicative since 
\item $f(p q) =f(p)f(q), \text{  }\gcd (p,q)=1$,
\item $f\left(p^2\right) \neq f(p)f(p), \text{     if }\ord_{p^2}(u)\neq p(p-1)$.
\end{enumerate}
\end{lem}

\begin{proof} The function has the value \(f\left(p^k\right)=1\) if and only if the element \(u\in \left(\mathbb{Z}\left/p^k\right.\mathbb{Z}\right)^{\times
	}\) is a primitive root modulo \(p^k\). Otherwise, it vanishes: \(f\left(p^k\right)=0\). The completely multiplicative property fails because of
	the existence of Wieferich primes, exempli gratia, \(0=f\left(40487^2\right) \neq f(40487)f(40487)=1\), see \cite{PA09}. Otherwise, it is completely multiplicative, that is, \(f\left(p^2\right) = f(p)f(p)=1\) for any nonWieferich primes $p \geq 2$. 
\end{proof}

Observe that the conditions \(\ord_p(u)=p-1 \text{ and } \ord_{p^2}(u)\neq p(p-1)\) imply that the integer \(u\neq \pm 1,v^2\) cannot be
extended to a primitive root \(\mod  p^k,k\geq 2\). But that the condition \(\ord_{p^2}(u)=p(p-1)\) implies that the integer \(u\) can
be extended to a primitive root \(\mod  p^k,k\geq 2\), see \cite[p.\ 208]{AT76}.

\section{Wirsing Formula} \label{s4}
This formula provides decompositions of some summatory multiplicative functions as products over the primes supports of the functions. This technique works well with certain multiplicative functions, which have supports on subsets of primes numbers of nonzero densities. \\

\begin{lem} {\normalfont (\cite[p. 71]{WE61})} \label{lem4.1}Suppose that \(f:\mathbb{N}\longrightarrow \mathbb{C}\) is a multiplicative function with the following properties.
\begin{enumerate} [font=\normalfont, label=(\roman*)]
\item $f(n) \geq 0$ for all integers $n\in \mathbb{N}$. 
\item $f\left(p^k\right)\leq c^k$ for all integers $k\in \mathbb{N}$, and $c<2$ constant. 
\item There is a constant $\tau >0$ such 
\begin{equation}
\sum _{p\leq x} f(p)=(\tau +o(1)) \frac{x}{\log  x}
\end{equation} 
as $x \longrightarrow  \infty$.  
\end{enumerate}

Then
\begin{equation}
	\sum _{n\leq x} f(n)=\left(\frac{1}{e^{\gamma \tau }\Gamma (\tau )}+o(1)\right)\frac{x}{\log  x}\prod _{p\leq x} \left(1+\frac{f(p)}{p}+\frac{f\left(p^2\right)}{p^2}+\cdots
	\right) .
	\end{equation}
	
\end{lem}

The gamma function appearing in the above formula is defined by \(\Gamma (s)=\int _0^{\infty }t^{s-1}e^{-s t}d t, s\in \mathbb{C}\). The intricate
proof of Wirsing formula appears in \cite{WE61}. It is also assembled in various papers, such as \cite{HA87}, \cite[p. 195]{PA88}, and discussed in \cite[p. 70]{MV07}, \cite[p. 308]{TG15}. Various applications are provided in \cite{MP11}, \cite{PS03}, \cite{WK75}, et alii. \\

\section{Harmonic Sums And Products Over Primes With Fixed Primitive Roots} \label{s5}
The subset of primes \(\mathcal{P}_u=\left\{ p\in \mathbb{P}:\text{ord}_p(u)=p-1 \right\}\subset \mathbb{P}\) consists of all the primes with a fixed primitive root \(u\in \mathbb{Z}\). By Hooley theorem, which is conditional on the generalized Riemann hypothesis, it has nonzero density \(\alpha
_u=\delta \left(\mathcal{P}_u\right)>0\). The real number \(\alpha _u>0\) coincides with the corresponding Artin constant, see \cite[p. 220]{HC67}, for the formula. The proof of the next result is based on standard analytic number theory methods in the literature, refer to \cite[Lemma 4]{PS03}.
\\

\begin{lem} \label{lem5.1} Assume the generalized Riemann hypothesis, and let \(x\geq 1\) be a large number. Then, there exists a pair
	of constants \(\beta _u>0,\text{     and     } \gamma _u>0\) such that
\begin{enumerate} [font=\normalfont, label=(\roman*)]
\item $ \displaystyle \sum_{\substack{p\leq x \\  p \in \mathcal{P}_u}} \frac{1}{p}=\alpha_u \log \log
	x+\beta_u+ O\left(\frac{\log \log x}{\log  x} \right) .$  
\item $ \displaystyle  
	\sum _{\substack{p\leq x \\ p\in \mathcal{P}_u}} \frac{\log p}{p-1}=\alpha_u \log
	x-\gamma_u+ O\left(\frac{\log \log x}{\log x}\right)$.
\end{enumerate}

\end{lem}

\begin{proof} (i). Let \(\pi _u(x)=\#\left\{ p\leq x:\text{ord}_p(u)=p-1 \right\}=\alpha _u\pi (x)\) be the counting measure of the
	corresponding subset of primes \(\mathcal{P}_u\). To estimate the asymptotic order of the prime harmonic sum, use the Stieltjes integral representation:
	\begin{equation}
		\sum _{\substack{p\leq x \\ p\in \mathcal{P}_u}} \frac{1}{p} =\int_{x_0}^{x}\frac{1}{t}d \pi _u(t) =\frac{\pi _u(x)}{x}+c_(x_0)
		+\int_{x_0}^{x}\frac{\pi _u(t)}{t^2}d t,
	\end{equation} 
	where \(x_0>0\) is a constant. Applying  Theorem \ref{thm7.1} yields
	\begin{eqnarray}
		\int_{x_0}^{x}\frac{1}{t}d \pi _u(t) 
		&=&\frac{\alpha _u}{\log  x}+O\left(\frac{\log \log x}{\log ^2(x)}\right)+c_0(x_0) \nonumber\\
		&\qquad&+\alpha
		_u\int_{x_0}^{x}\left(\frac{1}{t \log  t}+O\left(\frac{\log \log t}{t \log ^2(t)}\right)\right)d t \\
		&=&\alpha _u\log \log
		x-\log \log x_0+c_0\left(x_0\right)+O\left(\frac{\log \log x}{\log  x}\right)  \nonumber,
	\end{eqnarray} 
	where \(\beta _u=-\log \log x_0+c_0\left(x_0\right)\) is the Artin-Mertens constant. The statement (ii) follows from statement (i) and partial summation. 
\end{proof}

The Artin-Mertens constant \(\beta _u\) and the Artin-Euler constant \(\gamma _u\) have other equivalent definitions such as 
\begin{equation} \label{el330}
	\beta _u=\lim_{x \rightarrow
		\infty } \left(\sum _{p\leq x, p\in \mathcal{P}_u} \frac{1}{p}-\alpha _u\log \log x\right)\text{              }\text{   and   }\text{            }\beta
	_u=\gamma _u-\text{  }\sum _{p\in \mathcal{P}_u,} \sum _{k\geq 2} \frac{1}{k p^k} ,
\end{equation}
respectively. These constants satisfy \(\beta _u=\beta _1\alpha _u\text{ and } \gamma _u=\gamma \alpha _u\). If the density \(\alpha _u=1\), these definitions reduce 
to the usual Euler constant and the Mertens constant, which are defined by the limits 
\begin{equation} \gamma =\lim_{x \rightarrow  \infty }
	\left(\sum _{ p\leq x } \frac{\log  p}{p-1}-\log  x\right) \text{ and } \beta _1=\lim_{x \rightarrow  \infty } \left(\sum _{ p\leq x } \frac{1}{p}-\log \log x\right),
\end{equation}
or some other equivalent definitions, respectively. Moreover, the linear independence relation in (\ref{el330}) becomes \(\text{  }\beta =\gamma -\text{ 
}\sum _{p\geq 2} \sum _{k\geq 2} \left(k p^k\right) ^{-1}\), see \cite[Theorem 427] {HW08}.
\\

A numerical experiment for the primitive root \(u=2\) gives the approximate values 
\begin{enumerate} [font=\normalfont, label=(\roman*)]
\item $ \displaystyle  \alpha _2=\prod _{p\geq 2} \left(1-\frac{1}{p(p-1)}\right)=0.3739558667768911078453786 \ldots .$ 
\item $ \displaystyle \beta _2 \approx \sum_{ \substack{p\leq 1000 \\ p\in \mathcal{P}_2}} \frac{1}{p} -\alpha _2\log \log x=0.328644525584805374999956\text{ ... } , $  and 
\item  $ \displaystyle  \gamma _2 \approx \sum_{ \substack{p\leq 1000 \\ p\in \mathcal{P}_2}} \frac{\log  p}{p-1} -\alpha _2\log  x=0.424902273366234745796616 \ldots . $ 
\end{enumerate}

\begin{lem} \label{lem5.2} Assume the generalized Riemann hypothesis, and let \(x\geq 1\) be a large number. Then, there exists a pair of constants $\gamma _u>0$ and $\nu _u >0$ such that 
\begin{enumerate} [font=\normalfont, label=(\roman*)]
\item $\displaystyle  \prod _{\substack{ p\leq x\\ p\in \mathcal{P}_u}} \left(1-\frac{1}{p}\right)^{-1}=e^{\gamma
		_u} \log (x)^{\alpha _u}+O\left(\frac{\log \log x}{\log  x}\right) .$ 
\item $ \displaystyle \prod _{\substack{p\leq x\\ p\in \mathcal{P}_u}} \left(1+\frac{1}{p}\right)=e^{\gamma
		_u}\prod _{p\in \mathcal{P}_u} \left(1-p^{-2}\right) \log (x)^{\alpha _u}+O\left(\frac{\log \log x}{\log  x}\right) . $
\item $ \displaystyle \prod
	_{\substack{p\leq x\\ p\in \mathcal{P}_u}} \left(1-\frac{\log  p}{p-1}\right)^{-1}=e^{\nu _u-\gamma _u}x^{\alpha _u}+O\left(\frac{x^{\alpha _u}\log \log
		x}{\log  x}\right) .$ 
\end{enumerate}
\end{lem}

\begin{proof} (i). Express the logarithm of the product as 
	\begin{equation}
		\sum _{\substack{p\leq x\\ p\in \mathcal{P}_u}} \log \left(1-\frac{1}{p}\right)^{-1}=\sum
		_{\substack{p\leq x\\ p\in \mathcal{P}_u}} \sum _{k\geq 1} \frac{1}{k p^k}=\sum _{\substack{p\leq x\\ p\in \mathcal{P}_u}} \frac{1}{p} +\sum _{\substack{p\leq x\\ p\in \mathcal{P}_u}}
		\sum _{k\geq 2} \frac{1}{k p^k} .
	\end{equation}   
	
Apply Lemma \ref{lem5.1} to complete the verification. For statements (ii) and (iii), use similar methods as in the first one. 
\end{proof}

The constant \(\nu _u>0\) is defined by the double power series (an approximate numerical value for set \(\mathcal{P}_2=\{ 3, 5, 11, 13, \text{...}
\}\) is shown): 
\begin{equation}
	\nu _2=\sum _{p\in \mathcal{P}_2, } \sum _{k\geq 2} \frac{1}{k}\left(\frac{\log  p}{ p-1}\right)^k \approx 0.163507781570971567408003\text{... }.
\end{equation}\\

\section{Density Correction Factor} \label{s6}
The sporadic subsets of Abel-Wieferich primes, see \cite[p. 333]{RP96} for other details, have roles in the determination of the densities of the multiplicative subsets of integers \(\mathcal{N}_u\) with fixed primitive roots \(u\in \mathbb{Z}\), see Definition \ref{dfn1.05}. The prime product arising from the sporadic existence of the Abel-Wieferich primes \(p\geq 3\) is reformulated in the equivalent expression 
\begin{eqnarray} \label{el550}
	P(x)&=&\prod_{\substack{p^k\leq
			x,\\ \text{ord}_p(u)=p-1, \\ \text{ord}_{p^2}(u)\neq p(p-1)}} \left(1+\frac{1}{p}\right)\prod _{\substack{p^k\leq x,\\ \text{ord}_{p^2}(u)=p(p-1) }} \left(1+\frac{1}{p}+\frac{1}{p^2}+\cdots \right)  \nonumber\\
	&=&\prod _{\substack{p\leq x\\p\in \mathcal{W}}}  \left(1-\frac{1}{p^2}\right)\prod
	_{\substack{p\leq x\\p\in \mathcal{P}_u} } \left(1-\frac{1}{p}\right)^{-1} +O\left(\frac{1}{x}\right) \\
	&=&\prod _{p\in \mathcal{W}} \left(1-\frac{1}{p^2}\right)\prod _{\substack{p\leq x\\p\in \mathcal{P}_u} } \left(1-\frac{1}{p}\right)^{-1} +O\left(\frac{1}{x}\right)  \nonumber .
\end{eqnarray} 
Note that the subset of primes has the disjoint partition\\
\begin{equation}
	\mathcal{P}_u=\left\{ p\in \mathbb{P}:\text{ord}_p(u)=p-1 \right\}=\mathcal{W} \cup \overline{\mathcal{W}},
\end{equation}
where  
\begin{equation}
	\mathcal{W}=\left\{ p\in \mathbb{P}:\text{ord}_{p}(u)=p-1, \text{ ord}_{p^2}(u)\neq p(p-1)\right\}
\end{equation}
and
\begin{equation}
	\overline{\mathcal{W}}=\left\{ p\in \mathbb{P}:\text{ord}_{p^2}(u)=p(p-1)\right\}.
\end{equation}
The convergent partial product (\ref{el550}) is replaced with the approximation 
\begin{equation}
	\prod _{ \substack{p\leq x \\ p\in \mathcal{W}}} \left(1-\frac{1}{p^2}\right)=\prod
	_{p\in \mathcal{W}} \left(1-\frac{1}{p^2}\right) +O\left(\frac{1}{x}\right).
\end{equation} 

For \(u=2\), the subset of primes \(\mathcal{W}_2\) is the subset of Wieferich primes. This subset of primes is usually characterized in terms of the
congruence
\begin{equation}\label{eq6-50}
	\mathcal{W}_2=\left\{ p\in \mathbb{P}:2^{p-1}\equiv 1\text{ mod } p^2 \right\}=\{ 1093, 3511, \text{...}. \}. 
\end{equation}

Given a fixed \(u\neq \pm 1,v^2\), the product \(\prod _{ p\in \mathcal{W}} \left(1-p^{-2}\right)\) reduces the density to compensate for those primes
for which the primitive root \(u \text{ mod }p\) cannot be extended to a primitive root \(u \text{ mod }p^2\). This seems to be a density correction factor
similar to the case for primitive roots over the prime numbers. The correction factor required for certain densities of primes with respect to fixed primitive
roots over the primes was discovered by the Lehmers, see \cite{SP03}. \\

\section{The Proof Of The Theorem}  \label{s7}
The result below has served as the foundation for various other results about primitive roots. Most recently, it was used to prove the existence
of infinite sequences of primes with fixed prime roots, and bounded gaps, confer \cite{BP14}.\\

\begin{thm} { \normalfont (\cite{HC67})} \label{thm7.1} If it be assumed that the extended Riemann hypothesis hold for the Dedekind zeta function over Galois
	fields of the type \(\mathbb{Q}\left(\sqrt[d]{u},\sqrt[n]{1}\right)\), where \(n\) is a squarefree integer, and \(d \mid n.\) For a given nonzero integer
	\(u\neq \pm 1,v^2\), let \(A_u(x)\) be the number of primes \(p\leq x\) for which \(u\) is a primitive root modulo \(p\). Let \(u=u_1^m\cdot u_2^2\),
	where \(u_1>1\) is squarefree, and \(m\geq 1\) is odd. Then, there is a constant \(\alpha _u\geq 0\) such that
	\begin{equation}
		A_u(x)=\alpha _u\frac{x}{\log  x}+O\left(\frac{x
			\log \log x}{\log ^2 x}\right)\text{        }\text{    as } x\longrightarrow \infty  .
	\end{equation}
\end{thm}

\begin{proof} (Theorem \ref{thm1.2}) By the generalized Riemann hypothesis or Theorem \ref{thm7.1}, the density \(\alpha _u=\delta \left(\mathcal{P}_u\right)>0\) of the subset of primes \(\mathcal{P}_u\) is nonzero. Put \(\tau =\alpha _u\) in Wirsing formula, Lemma \ref{lem4.1}, and replace the characteristic function
\(f(n)\) of primitive roots in the finite ring \(\mathbb{Z}\left/p^k\right. \mathbb{Z},k\geq 1\), see Lemma \ref{lem4.1}, to produce\\
\begin{eqnarray} \label{700}
	N_u(x)&=&\sum _{n\leq x} f(n)\nonumber \\
&=&\left(\frac{1}{e^{\gamma \tau }\Gamma (\tau )}+o(1)\right)\frac{x}{\log x}\prod _{p^k\leq x} \left(1+\frac{f(p)}{p}+\frac{f\left(p^2\right)}{p^2}+\cdots \right)   \nonumber\\
	&=&\left(\frac{1}{e^{\gamma \alpha _u}\Gamma \left(\alpha _u\right)}+o(1)\right)\frac{x}{\log  x} \\
	&\qquad& \times \prod _{\substack{ p^k\leq x, \\ \text{ord}(u)=p-1,\\ \text{ord}(u)\neq
			p(p-1)} } \left(1+\frac{1}{p}\right)\prod _{\substack{p^k\leq x, \\ \text{ord}(u)=p(p-1)}} \left(1+\frac{1}{p}+\frac{1}{p^2}+\cdots \right)   \nonumber.
\end{eqnarray}
In equation (\ref{700}), line 1, the product over the prime powers $p^k \leq x$ is broken up into two subproducts. In line 2, the first subproduct is restricted to the subset of Wieferich prime powers which do not satisfy the completely multiplicative property $f(p^2)\ne f(p)f(p)$ of the characteristic function; the second subproduct is restricted to the subset of nonWieferich prime powers which do satisfy the completely multiplicative property $f(p^2)= f(p)f(p)$ of the characteristic function, Lemma \ref{lem3.1}.\\

Replacing the equivalent product, see (\ref{el550}) in Section \ref{s6}, and using Lemma \ref{lem5.2}, yield \\
\begin{eqnarray}
	\sum _{n\leq x} f(n)&=&\left(\frac{1}{e^{\gamma \alpha _u}\Gamma
		\left(\alpha _u\right)}+o(1)\right)\frac{x}{\log  x}\prod _{ p\in \mathcal{W}} \left(1-\frac{1}{p^2}\right)\prod _{\substack{p^k\leq x\\ p\in \mathcal{P}_u}}
	\left(1-\frac{1}{p}\right) ^{-1}   \nonumber\\
	&=&\left(\frac{e^{\gamma _u-\gamma \alpha _u}}{\Gamma \left(\alpha _u\right)}+o(1)\right)\frac{x}{(\log  x)^{1-\alpha
			_u}}\prod _{p\in \mathcal{W}} \left(1-\frac{1}{p^2}\right)  ,
\end{eqnarray}\\
where \(\gamma _u\) is the Artin-Euler constant, see Lemmas \ref{lem5.1} and \ref{lem5.2} for details. Lastly, the error term \(o\left(x(\log  x)^{\alpha_u-1}\right)\)
absorbs all the errors.          Quod erat demonstrandum. 
\end{proof}

\section{Harmonic Sum For The Fixed Primitive Root 2} \label{s8} 
Let \(\mathcal{N}_2\) be the subset of integers such that \(u=2\) is a primitive root
modulo \(n\geq 1\), and let \(N_2(x)=\#\left\{ n\leq x:n\in \mathcal{N}_2 \right\}\) be the corresponding discrete counting measure, see Theorem \ref{thm1.1}. 
The subset \(\mathcal{N}_2\subset \mathcal{A}\) is a proper subset of \(\mathcal{A}=\{ n\in \mathbb{N}:p \mid n\Rightarrow p\in \mathcal{P}_2\)$\}$, which is 
generated by the subset of primes \(\mathcal{P}_2=\left\{ p\in \mathbb{P}:\ord_p(2)=p-1 \right\}\). Since 2 is not a primitive root modulo the prime powers \(p_0^m,m\geq 2\), the subset \(\mathcal{A}\) is slightly larger than the subset \(\mathcal{N}_2\). More precisely, the Wieferich prime powers \(p_0^m\in \mathcal{A},m\geq 1\), but \(p_0^m\notin \mathcal{N}_2,m\geq 2\), where $p_0^{p-1}-1 \equiv 0 \bmod p_0^2$. \\
 
An asymptotic formula for the harmonic sum over the subset of integers \(\mathcal{N}_2\) is determined here.

\begin{lem} \label{lem9.1}  Let \(x\geq 1\) be a large number, let \(\alpha _2=\delta \left(\mathcal{P}_2\right)>0\) be the density of the subset of primes
 \(\mathcal{P}_2\), and let \(\mathcal{N}_2\subset \mathbb{N}\) be a subset of integers generated by \(\mathcal{P}_2\). Then
\begin{equation}
\sum _{\substack{n \leq  x\\ n\in \mathcal{N}_2}} \frac{1}{n}=\kappa _2 (\log  x)^{\alpha _2}+\gamma _2+O\left(\frac{1}{(\log  x)^{1-\alpha _2}}\right) .
\end{equation} 
The number \(\alpha _2=0.373955\text{...}\) is Artin constant, and \(\gamma _2>0\) is the Artin-Euler constant, see (\ref{el330}) for the definition. The other constant is 
\begin{equation}
\kappa _2=\frac{e^{\gamma _2-\gamma \alpha _2}}{\alpha _2\Gamma \left(\alpha _2\right)}\prod _{p\in \mathcal{W}} \left(1-\frac{1}{p^2}\right) \approx 1.12486444988498798741328\text{...}.
\end{equation} 
This numerical approximation assume that $\gamma _2-\gamma \alpha _2=0$, and the index of the product ranges over the subset Wieferich primes \(\mathcal{W}\). 
\end{lem}

\begin{proof} Use the discrete counting measure \(N_2(x)=\left(\alpha _2\kappa _2+o(1)\right)x(\log  x)^{\alpha _2-1}\), Theorem 1.1, to write the finite sum as an integral, and evaluate it: 
\begin{equation}
\left .\sum _{\substack{n \leq  x\\ n\in \mathcal{N}_2}} \frac{1}{n} =\int_{x_0}^{x}\frac{1}{t}d N_2(t) =\frac{N_2(t)}{t} \right | _{x_0}^{x}+\int_{x_0}^{x}\frac{N_2(t)}{t^2}d t,
\end{equation} 
where \(x_0>0\) is a constant. Continuing the evaluation yields 
\begin{eqnarray}
\int_{x_0}^{x}\frac{1}{t}d N_2(t) 
&=&
\frac{\left(\alpha _2\kappa
_2+o(1)\right)}{(\log  x)^{1-\alpha _2}}+c_0\left(x_0\right)+\int_{x_0}^{x} \frac{\left(\alpha _2\kappa _2+o(1)\right)}{t(\log
 t)^{1-\alpha _2}}d t \\
&=&\kappa _2(\log  x)^{\alpha _2}+\gamma _2+O\left(\frac{1}{(\log  x)^{1-\alpha _2}}\right)  \nonumber,
\end{eqnarray}
where \(c_0\left(x_0\right)\) is a constant. Moreover, 
\begin{equation}
\gamma _2=\lim_{x \rightarrow  \infty } \left(\sum _{\substack{n \leq  x\\ n\in \mathcal{N}_2}} \frac{1}{n}-\kappa
_2\log ^{\alpha _2} x\right)=c_0\left(x_0\right)+\int_{x_0}^{\infty}\frac{\left(\alpha _2\kappa _2+o(1)\right)}{t(\log  t)^{1-\alpha
_2}}d t
\end{equation} 
is a second definition of this constant. 
\end{proof}

The integral lower limit \(x_0=2\) appears to be correct one since the subset of integers is \(\mathcal{N}_2=\left\{ 3, 5, 3^2,11,\text{...} \right\}\).



\begin{thebibliography}{999}
\bibitem{PA09} Paszkiewicz, A. A new prime $p$ for which the least primitive root $\bmod p$ and the least primitive root $\bmod p^2$   are not equal. Math. Comp.  78  (2009),  no. 266, 1193-1195.
\end{thebibliography}
\end{document}